\documentclass[11pt]{amsart}
\usepackage{amsmath,amssymb,amsfonts,amsthm,amscd,indentfirst}

\usepackage[dvipdfm]{hyperref}

\newtheorem{prop}{Proposition}
\newtheorem{theo}[prop]{Theorem}
\newtheorem{lemm}[prop]{Lemma}

\newtheorem{defi}[prop]{Definition}

\def\begeq{\begin{equation}}
\def\endeq{\end{equation}}

\def\and{\quad{\rm and}\quad}

\def\<{\langle}
\def\>{\rangle}

\def\Dint{\displaystyle\int}
\def\Dfrac{\displaystyle\frac}

\begin{document}
\title[The quermassintegral inequalities]{The quermassintegral
inequalities for $k$-convex starshaped domains}
\author{Pengfei Guan }
\address{Department of Mathematics\\
         McGill University\\
         Montreal, Quebec. H3A 2K6, Canada.}
\email{guan@math.mcgill.ca}
\author{Junfang Li}
\address{Department of Mathematics\\
         McGill University\\
         Montreal, Quebec. H3A 2K6, Canada.}
\email{jli@math.mcgill.ca}
\thanks{Research of the first author was supported in part by an
NSERC Discovery Grant, and research of the second author is
supported in part by a CRM fellowship.}
\begin{abstract}
We give a proof of the isoperimetric inequality for quermassintegrals of non-convex
starshaped domains, using a result of Gerhardt \cite{G} and Urbas \cite{U}
on an expanding geometric curvature flow.
\end{abstract}

\date{}
\maketitle

The Alexandrov-Fenchel inequalities \cite{Alex37, Alex38} for the
quermassintegrals of convex domains are fundamental in classical
geometry. For a bounded domain $\Omega\subset \mathbb R^{n+1}$, we
denote $M=\partial \Omega$ the boundary of $\Omega$. We will assume
$M$ smooth in this paper. Let
\[\kappa(x)=(\kappa_1(x),\cdots, \kappa_n(x))\] be the principal curvatures of $x\in
M$, and let $\sigma_k(\lambda)$ the $k$th elementary function in
$\lambda=(\lambda_1,\cdots, \lambda_n) \in \mathbb R^n$ (with
$\sigma_0(\lambda)\equiv 1$). There are several equivalent definitions of the
quermassintegral $V_{(n+1)-k}(\Omega)$. For positive integer $k$, we will take the following
\begin{equation}\label{quer-def0}
V_{(n+1)-k}(\Omega)=C_{n,k}\int_{M}
\sigma_{k-1}(\kappa)d\mu_{M},\end{equation} where
where $\sigma_k$ is the
$k$th elementary symmetric function,
\begin{equation}\label{Cnk}
C_{n,k}=\frac{\sigma_k(I)}{\sigma_{k-1}(I)},\end{equation} with
$I=(1,\cdots,1)$. One may also recover $V_{n+1}(\Omega)$ by the Minkowski type formula,
\begin{equation}\label{quer-def}
V_{(n+1)-k}(\Omega)=\int_{M}
u\sigma_k(\kappa)d\mu_{M},\end{equation}
where $u=\langle X,\nu\rangle$, $X$
is the position function of $M$, and $\nu$ is the outer-normal of
$M$ at $X$.
$V_{n+1}(\Omega)$ is a multiple of the volume of
$\Omega$ by a dimensional constant, $V_n(\Omega)$ is a multiple of
the surface area of $M=\partial \Omega$ by another dimensional
constant.
If $\Omega$ is convex, the celebrated Alexandrov-Fenchel
quermassintegral inequality states that, for $0\le m\le n$,
\begin{equation}\label{AFI}
\Big(\frac{V_{(n+1)-m}(\Omega)}{V_{(n+1)-m}(B)}\Big)^{\frac{1}{n+1-m}}\le
\Big(\frac{V_{n-m}(\Omega)}{V_{n-m}(B)}\Big)^{\frac{1}{n-m}},
\end{equation} where $B$ is the standard ball in $\mathbb R^{n+1}$.
The equality holds if and only if $\Omega$ is a ball. The case $m=0$
is the classical isoperimetric inequality.

There have been some interests in extending the original
Alexandrov-Fenchel inequality to non-convex domains (e.g.,
\cite{Tru}, \cite{Gibbons}). In this short paper, we extend this
inequality to starshaped domains in $\mathbb R^{n+1}$. These domains
are special type of domains where fully nonlinear partial
differential equations were studied in pioneer work by
Caffarelli-Nirenberg-Spruck \cite{CNS1, CNS}. We follow the
notations in \cite{Tru} as below.
\begin{defi}For $\Omega\subset \mathbb R^{n+1}$, we say $\Omega$ is
$k$-convex if $\kappa(x)\in \bar \Gamma_k$
for all $x\in M$, where $\Gamma_k$ is the Garding's cone
\[\Gamma_k=\{\lambda \in \mathbb R^n \ | \quad \sigma_m(\lambda)>0, \quad  \forall m\le k\}.\]
We say $\Omega$ is strictly $k$-convex if $\kappa(x)\in \Gamma_k$
for all $x\in M$. \end{defi}
$n$-convex is {\it convex} in usual sense, $1$-convex is sometimes
referred as {\it mean convex}.
\begin{theo}\label{tAFI} Suppose $\Omega$ is a $k$-convex starshaped domain in
$\mathbb R^{n+1}$, then inequality (\ref{AFI}) is true for $0\le
m\le k$. The equality holds if and only if $\Omega$ is a
ball.\end{theo}

In \cite{Tru}, Trudinger considered inequality (\ref{AFI}) for
general the $k$-convex domains. He proposed an elliptic method by
reducing the problem to a Hessian type equation in the domain, but
the reduction argument there is incomplete. Our proof here is a
parabolic one, using the flow studied by Gerhardt \cite{G} and Urbas
\cite{U}. For our purpose, we will only use a special case of their
result for the following evolution equation on a hypersurface $M_0$
in $\mathbb{R}^{n+1 }$,
\begin{equation}\label{GGGe}
X_t=\frac{\sigma_{k-1}}{\sigma_{k}} (\kappa) \nu.\end{equation}
\begin{theo}\label{tG-U} (Gerhardt \cite{G}, Urbas \cite{U})
If $\Omega_0$ is a starshaped strictly $k$-convex domain, then solution for
flow (\ref{GGGe}) exists for all time $t>0$ and it converges to a
round sphere after a proper rescaling.\end{theo}

\medskip
We define
\begin{equation}\label{iso-ratio} \mathcal {I}_k(\Omega)
=\frac{V_{(n+1)-k}^{\frac{1}{n+1-k}}(\Omega)}{V_{n-k}^{\frac{1}{n-k}}(\Omega)}.\end{equation}
The key observation is that the isoperimetric ratio $\mathcal
{I}_k(\Omega)$ of the quermassintegrals are monotone along expanding
flow of (\ref{GGGe}). From what follows, we will denote $M(t)$ the
solution of flow (\ref{GGGe}) at time $t$. If there is no confusion,
we will just write $M(t)=M$. To simplify notation, we will also
write $\sigma_m$ for $\sigma_m(\kappa)$ unless specified otherwise.
To prepare our proof of Theorem \ref{tAFI}, we first list the
evolution equations of various geometric quantities under the
following general evolution equation.
\begin{equation}\label{GEE gen}
\begin{array}{rll}
\partial_t X=F\nu,
\end{array}
\end{equation}
where $F=F(\kappa, X, t)$.
\begin{prop}\label{prop1-1}
Under flow (\ref{GEE gen}),  we have the following evolution equations.
\begin{equation}\label{coro 1 sect 2 equ1}
\begin{array}{rll}
&&\partial_t g_{ij} = 2Fh_{ij}\\
&&\partial_t\nu = -\nabla F\\
&&\partial_t (d\mu_g)  = F\sigma_1 d\mu_g\\
&&\partial_t h_{ij} = -\nabla_i\nabla_j F + F (h^2)_{ij}\\
&&\partial_t h^i_j  =
-\nabla^i\nabla_j F - F (h^2)^i_j\\
&&\partial_t \sigma_m = -\nabla_j([T_{m-1}]^i_j\nabla_i F)
- F\sigma_{m-1,1}(h^i_j;(h^2)^i_j)\\
\end{array}
\end{equation}
where we denote $h^i_j\equiv g^{ik}h_{kj}$, $(h^2)_{ij}\equiv
g^{kl}h_{ik}h_{lj}$ and $(h^2)^i_{j}\equiv
g^{is}g^{kl}h_{sk}h_{lj}$, $[T_{l}]^i_j$ is the $l$-th Newton
transformation, and $\sigma_{m-1,1}(A;B)=\frac{\partial
\sigma_m}{\partial A_{ij}}(A)B_{ij}$ is a polarization of
$\sigma_m$.
\end{prop}

\begin{proof}
The proof follows from straightforward computations using the
Codazzi property of the second fundamental form. The last identity
follows from the divergent free property of $[T_{m-1}]^i_j$.
\end{proof}

\begin{lemm}\label{coro 1 sect 2}
Under flow (\ref{GGGe}),
\begin{equation}\label{1.3}
\partial_t \Dint_M \sigma_l d\mu_g = (l+1)\Dint_M \Dfrac{\sigma_{l+1}\sigma_{k-1}}{\sigma_k}
d\mu_g.
\end{equation}
\end{lemm}

\begin{proof}
From identities in Proposition \ref{prop1-1}, for $1\le l\le n$, we
have
\begin{equation}\label{1.3-1}
\begin{array}{rll}
\partial_t \Dint_M \sigma_l d\mu_g &=
\Dint_M \partial_t\sigma_ld\mu_{g} +\sigma_l \partial_t
(d\mu_{g})\\
&= -\Dint_M \frac{\sigma_{k-1}}{\sigma_k}\bigg( \sigma_{l-1,1}(h^i_j;(h^2)^i_j)-\sigma_l\sigma_1 \bigg)d\mu_g\\
&= (l+1)\Dint_M \frac{\sigma_{k-1}}{\sigma_k}\sigma_{l+1}d\mu_g,\\
\end{array}
\end{equation}
where we have used identity \begin{equation}
\sigma_{l-1,1}(h^i_j;(h^2)^i_j)=\sigma_1\sigma_l-(l+1)\sigma_{l+1}.\end{equation}\end{proof}

\medskip

There exist literatures using flow to establish geometric
inequalities (e.g. \cite{Andrews, GW1, GW2}). To compare two
geometric quantities, one would like to design a normalized flow so
that one of the quantities is invariant under the flow, and another
is monotone along the flow. Suppose we consider flow (\ref{GGGe}),
and throw a time dependent constant $R(t)$ to normalize it as the
following
\begin{equation}\label{GGGe-N}
 X_t=(\frac{\sigma_{k-1}}{\sigma_{k}}(\kappa) -R(t)) \nu.\end{equation}
Suppose $\Omega_t$ is the domain with $X(t)$ as position function.
The normalization constant $R(t)$ should be picked to make
$V_{n-k}(\Omega_t)\equiv constant$ along (\ref{GGGe-N}), which can
be calculated easily using previous lemma. But one runs in to
trouble to establish a priori estimates for flow (\ref{GGGe-N})
under $k$-convexity assumption. By invoking the Minkowski identity,
the right normalized flow should be
\begin{equation}\label{GGGe-R}
 X_t=(\frac{\sigma_{k-1}}{\sigma_{k}}(\kappa) -r(t)u) \nu,\end{equation}
where $u$ is the support function of $M(t)$, and
\begin{equation}\label{r(t)} r(t)=\frac{\int_{M} \frac{\sigma_{k+1}\sigma_{k-1}}{\sigma_k}
d\mu_g}{C_{n,k+1}\int_M \sigma_{k}d\mu_g}.\end{equation}
It is straightforward to show that $r(t)$ is a
normalization constant to make $V_{n-k}(\Omega_t)$ invariant
under the flow, and $V_{n-k+1}(\Omega_t)$ is nondecreasing! Plus, one may establish all
the a priori estimates for the normalized flow (\ref{GGGe-R}). Inequality (\ref{AFI}) can be proved
along the way.

On the other hand, it turns out that flow (\ref{GGGe-R}) is equivalent (up to an isomorphism) to
\begin{equation}
 X_t=\frac{\sigma_{k-1}}{\sigma_{k}}(\kappa) \nu -r(t)X,\end{equation}
which in turn is a re-parametrization of the original flow (\ref{GGGe}). Therefore,
we have the following simple proof using directly the result of Gerhardt and Urbas in Theorem \ref{tG-U}.

\medskip

\noindent{\bf Proof of Theorem \ref{tAFI}}.
It is easy to see that $\mathcal {I}_k(\Omega)$ is invariant under
rescaling. We only need to show that,
\begin{equation}\label{iso-inc} \mathcal {I}_k(\Omega)\le \mathcal
{I}_k(B),\end{equation} and the equality holds if and only if
$\Omega$ is a ball.

\noindent
{\it Case 1. $\Omega$ is strictly $k$-convex.}

For solution $X(\cdot,t)$ in Theorem \ref{tG-U}, consider $\tilde
X(\cdot,t)=e^{-\int_0^tr(s)ds}X(\cdot,t)$, where $r(t)$ as in (\ref{r(t)}).
We denote
$\tilde \Omega_t$ to be the domain enclosed by $\tilde X(\cdot,t)$.

Since $X(\cdot,t)$ is converging to a sphere (after a proper
rescaling), we only need to show $\mathcal {I}_k(\tilde \Omega_t)$
is increasing. We will continue to denote
$\sigma_m=\sigma_m(\kappa)$, where $\kappa$ is the principal
curvature of $X$.

From (\ref{1.3}) in Lemma \ref{coro 1 sect 2}, with $C_{n,k}$
defined as in (\ref{Cnk}), we have
\begin{equation}\begin{array}{rll}
\frac{d V_{n-k}(\tilde \Omega_t)}{dt} &=
(k+1)C_{n,k+1}e^{-(n-k)\int_0^t r(s)ds}\bigg[\Dint_M
\frac{\sigma_{k-1}}{\sigma_k}
\sigma_{k+1}d\mu_g-rC_{n,k+1}\Dint_M \sigma_{k}d\mu_g\bigg]\\
&= 0,
\end{array}
\end{equation} and
\begin{equation}\label{00}
\begin{array}{rll}
\frac{d V_{(n+1)-k}(\tilde \Omega_t)}{dt} &=
kC_{n,k}e^{-(n+1-k)\int_0^t r(s)ds}\bigg[\Dint_M
\frac{\sigma_{k-1}}{\sigma_k}
\sigma_{k}d\mu_g-rC_{n,k}\Dint_M \sigma_{k-1}d\mu_g\bigg]\\
&= kC_{n,k}e^{-(n+1-k)\int_0^t r(s)ds}\Dint_M \bigg[1-\frac{\int_M
\frac{\sigma_{k+1}\sigma_{k-1}}{\sigma_k}
d\mu_g}{C_{n,k+1}\int_M \sigma_{k}d\mu_g}C_{n,k} \bigg]\sigma_{k-1}d\mu_g \\
&\ge kC_{n,k}e^{-(n+1-k)\int_0^t r(s)ds}\Dint_M
\bigg[1-\Dfrac{\sigma_{k+1}(I)\sigma_{k-1}(I)}{\sigma^2_k(I)}
\Dfrac{C_{n,k}}{C_{n,k+1}}
\bigg]\sigma_{k-1}d\mu_g\\
&=0,\\
\end{array}
\end{equation}
where we have used the Newton-MacLaurin inequality in the last step
of (\ref{00}). If the equality holds in (\ref{iso-inc}), we must
have $\frac{d V_{(n+1)-k}(\tilde \Omega_t)}{dt}\equiv 0$. Therefore,
the equality of the Newton-MacLaurin inequality must be held at
every point of $M$ in (\ref{00}). This implies $M$ is a round sphere
for each $t\ge 0$. In particular, $M_0$ is a sphere.

\noindent
{\it Case 2. General $k$-convex starshaped domain $\Omega$.}

We may approximate it by strictly $k$-convex starshaped domains. The
inequality follows from the approximation. We now treat the equality
case. We first note that both $\int_M\sigma_kd\mu_g$ and
$\int_M\sigma_{k-1}d\mu_g$ are positive, since there exists at least
one elliptic point on an embedded compact hypersurface in Euclidean
space and also the $k$-convexity condition. Suppose $\Omega$ is a
$k$-convex starshaped domain with equality in (\ref{iso-inc})
attained. Let $M_{+}=\{x\in M | \sigma_k(\kappa(x))>0\}$. $M_{+}$ is
open and nonempty since $M$ is compact and embedded in $\mathbb
R^{n+1}$. We claim that $M_{+}$ is closed. This would imply
$M=M_{+}$, so $\Omega$ is strictly $k$-convex, by {\it Case 1}, we
may conclude $\Omega$ is a standard ball.

We now prove that $M_{+}$ is closed. Pick any $\rho\in C^2_0(M_{+})$
compactly supported in $M_{+}$. Let $M_s$ be the hypersurface
determined by position function $X_s=X+s\rho \nu$, where $X$ is the
support function of $M$ and $\nu$ is the unit outernormal of $M$ at
$X$. Let $\Omega_s$ be the domain enclosed by $M_s$. It is easy to
show $M_s$ is $k$-convex starshaped when $s$ is small enough.
Therefore $\mathcal {I}_k(\Omega_s)-\mathcal {I}_k(\Omega)\le 0$ for
$s$ small, i.e.
$$\frac{d}{ds}\mathcal {I}_k(\Omega_s)|_{s=0}= 0.$$
Simple calculation yields
\[\frac{d}{ds} \Dint_{M_s}\sigma_l(\kappa_s) d\mu_{g_s} |_{s=0}=
(l+1)\Dint_M \sigma_{l+1}(\kappa)\rho d\mu_g.\]
Therefore, we have
\[\frac{d}{ds}\mathcal {I}_k(\Omega_s)|_{s=0}
=A\Dint_M (\sigma_{k+1}(\kappa)-c_1 \sigma_k(\kappa))\rho d\mu_g= 0,\]
for some constant $A>0$ with $c_1
=\frac{k(n-k)}{(k+1)(n-k+1)}\frac{1}{\mathcal{I}(B)^{n-k+1}
(\int_M\sigma_k)^\frac{1}{n-k}}>0$ and for all $\rho\in
C^2_0(M_{+})$.

In turn,
\begin{equation}\label{qq1} \sigma_{k+1}(\kappa(x))=c_1 \sigma_k(\kappa(x)),
\quad \quad \forall x\in M_{+}.\end{equation}
By the Newton-MacLaurine inequality, there is a dimensional constant
$\tilde C_{k,n}$ such that $\sigma_{k+1}(\kappa(x))\le \tilde C_{k,n}\sigma^{1+1/k}_{k}(\kappa(x))$
for all $x\in M_{+}$. In view of (\ref{qq1}), there is a positive constant $c_2$, such that
\begin{equation}\label{qq2} \sigma_k(\kappa(x))\ge c_2>0, \quad \quad \forall x\in M_{+},\end{equation}
where $c_2=(\frac{c_1}{\tilde C_{k,n}})^k$ is a positive constant
depending only on $n$, $k$, and $\Omega$. (\ref{qq2}) implies
$M_{+}$ is closed. \qed

\medskip

The question of validity of inequality (\ref{AFI}) for general $k$-convex domains is
still open. When $k=1$, flow (\ref{GGGe}) is exactly the inverse
mean curvature flow. There is a notion of weak solution studied by Huisken-Ilmanen
for the Penrose inequality in \cite{HI}. Huisken \cite{Hu}
has proved the following with outward minimising assumption.

\begin{theo}\label{tHu} If $\Omega \subset R^{n+1}$ is outward minimising for
$n\leq 6$, then inequality (\ref{AFI}) is valid for $m=1$. \end{theo}

\noindent {\bf Acknowledgment:} We would like to
thank Bennett Chow and Guofang Wang for several enlightening discussions. We would also
like to thank Gerhard Huisken for informing us his result in Theorem \ref{tHu}.

\end{document}